\documentclass[11pt, reqno]{amsart}
\usepackage{amsmath,amsthm,amssymb}

\theoremstyle{plain}
\newtheorem{thm}{Theorem}[section]
\newtheorem{prop}[thm]{Proposition}
\newtheorem{lem}[thm]{Lemma}
\newtheorem{cor}[thm]{Corollary}

\theoremstyle{definition}

\newtheorem{rem}[thm]{Remark}
\newtheorem{defn}[thm]{Definition}
\newtheorem{eg}[thm]{Example}
\newtheorem{subtitle}[thm]{}
\newtheorem{ex}{Exercise}[section]
\numberwithin{equation}{section}

\def\b{\beta}

\def\K{\nabla}
\def\l{\lambda}

\def\n{\,\vert\,}
\def\o{\theta}
\def\w{\omega}

\def\cg{{\mathcal{G}}}

\def\cl{{\mathcal{L}}}

\def\co{{\mathcal{O}}}

\def\li{\langle}
\def\ri{\rangle}
\def\n{\ \vert\ }

\def\ms{\medskip}
\def\ss{\smallskip}

\def\ni{\noindent}
\def\ti{\tilde}
\def\p{\partial}

\def\I{{\rm I\/}}
\def\II{{\rm II\/}}

\def\R{\mathbb{R} }
\def\C{\mathbb{C}}

\def\Z{\mathbb{Z}}
\def\O{\mathbb{O}}
\def\fg{\mathfrak{G}}

\newcommand{\beg}{\begin{eg}}
\newcommand{\eeg}{\end{eg}}
\newcommand{\bthm}{\begin{thm}}
\newcommand{\ethm}{\end{thm}}
\newcommand{\bprop}{\begin{prop}}
\newcommand{\eprop}{\end{prop}}
\newcommand{\bcor}{\begin{cor}}
\newcommand{\ecor}{\end{cor}}
\newcommand{\blem}{\begin{lem}}
\newcommand{\elem}{\end{lem}}
\newcommand{\bca}{\begin{cases}}
\newcommand{\eca}{\end{cases}}
\newcommand{\brem}{\begin{rem}}
\newcommand{\erem}{\end{rem}}
\newcommand{\bpm}{\begin{pmatrix}}
\newcommand{\epm}{\end{pmatrix}}
\newcommand{\bbm}{\begin{bmatrix}}
\newcommand{\ebm}{\end{bmatrix}}
\newcommand{\bvm}{\begin{vmatrix}}
\newcommand{\evm}{\end{vmatrix}}
\newcommand{\bdefn}{\begin{defn}}
\newcommand{\edefn}{\end{defn}}
\newcommand{\bsub}{\begin{subtitle}}
\newcommand{\esub}{\end{subtitle}}
\newcommand{\bex}{\begin{ex}}
\newcommand{\eex}{\end{ex}}
\newcommand{\ben}{\begin{enumerate}}
\newcommand{\een}{\end{enumerate}}

\newcommand{\beq}{\begin{equation}}
\newcommand{\eeq}{\end{equation}}
\newcommand{\beqn}{\begin{eqnarray*}}
\newcommand{\eeqn}{\end{eqnarray*}}
\newcommand{\beqa}{\begin{eqnarray}}
\newcommand{\eeqa}{\end{eqnarray}}
\newcommand{\CR}{\nonumber \\}

\def\rmG{\mathrm{G}}
\def\rmS{\mathrm{S}}
\def\rmT{\mathrm{T}}

\def\rmSO{\mathrm{SO}}

\def\rmd{\mathrm{d}}
\def\rmG{\mathrm{G}}
\def\rmT{\mathrm{T}}
\def\rmspan{\mathrm{span}}
\def\fg{\frak{g}}
\def\fh{\frak{h}}
\def\fb{\frak{b}}

\author{Shengli Kong}
\address{Department of Mathematics\\
University of California at Irvine, Irvine, CA 92697-3875}
\email{skong@math.uci.edu}
\author{Chuu-Lian Terng$^\ast$}\thanks{$^\ast$Research supported
in  part by NSF grant DMS-0529756}
\address{Department of Mathematics\\
University of California at Irvine, Irvine, CA 92697-3875}
\email{cterng@math.uci.edu}
\author{Erxiao Wang}
\address{Department of Mathematics\\
  University of Texas at Austin, Austin, TX 78712-0257}
\email{ewang@math.utexas.edu}

\begin{document}

\title[Associative Cones and Integrable Systems]
{Associative Cones and Integrable Systems}


\maketitle

\centerline{\it Dedicated to the memory of Shing-Shin Chern}

\begin{abstract}
We identify $\R^7$ as the  pure imaginary part of octonions. Then the multiplication in octonions gives a natural almost complex structure for the unit sphere $S^6$.  It is known that a cone over a surface $M$ in $\rmS^6$ is an associative submanifold of $\R^7$ if and only if $M$ is almost complex in $\rmS^6$.  In this paper, we show that the Gauss-Codazzi equation for almost complex curves in $S^6$ is the equation for primitive maps associated to the
$6$-symmetric space $G_2/T^2$, and use this to explain some of the
known results.  Moreover, the equation for $\rmS^1$-symmetric almost
complex curves in $S^6$ is the periodic Toda lattice associated to $G_2$, and a
discussion of periodic solutions is given. 

\end{abstract}

\section{Introduction}

We identify $\R^7$ as the pure imaginary part of the octonions $\O$.
It is known that the group of automorphism of $\O$ is the compact
simple Lie group $G_2$, and the constant $3$-form  on $\R^7$,
$$\phi(u_1, u_2, u_3)= (u_1\cdot u_2, u_3),$$
is invariant under $G_2$.     A $3$-dimensional submanifold $M$ in
$\R^7$ is {\it associative\/} if $\R 1+ TM_x$ is an associative
subalgebra of $\O$ for all $x\in M$, i.e., it is isomorphic to the
quaternions.  It is easy to see that a $3$-dimensional submanifold
of $\R^7$ is associative if and only if it is calibrated by the
$3$-form $\phi$.

The multiplication of octonions defines an almost complex structure
on the unit sphere $S^6$ by $J_x(v)= x\cdot v$.  An immersion $f$
from a Riemann surface $\Sigma$ to $S^6$ is called {\it almost
complex\/} if the differential of $f$ is complex linear, i.e.,
$df_x(iv)= J_x(df_x(v))= x\cdot df_x(v)$. It is known that
(\cite{HL}) a surface $\Sigma$ is an almost complex curve in $S^6$
if and only if the cone over $\Sigma$ is an associative submanifold
of $\R^7$.

   An immersion $f$ from a Riemann surface to $S^n$ is called
   {\it totally isotropic\/} if $((\K_{\frac{\p}{\p
z}})^if_\ast(\frac{\p}{\p z}), (\K_{\frac{\p}{\p
z}})^jf_\ast(\frac{\p}{\p z}))=0$
   for all $i, j\geq 0$, where $(X, Y)=\sum_{i=1}^{n+1} X_i Y_i$ is the
complex bilinear form on $\C^{n+1}$.   A surface in $S^n$ is said to
be {\it full}
   if it does not contain in any hypersphere.
Bolton, Vrancken, and Woodward (\cite{BVW}) used harmonic sequences
to prove that if $f:\Sigma\to S^6$ is an immersed almost complex
curve, then $f$ must be one of the following:
\begin{itemize}
   \item[(i)] full in $\rmS^6$ and totally isotropic,
   \item[(ii)] full in $\rmS^6$ and not totally isotropic,
   \item[(iii)] full in some totally geodesic $\rmS^5$ in $\rmS^6$,
   \item[(iv)] a totally geodesic $S^2$.
\end{itemize}
Bryant (\cite{B0}) used twistor theory to construct type (i) almost
complex curves of any genus in $S^6$. Cones over a type (iii) almost
complex curves in $S^6$ are  special Lagrangian submanifolds, which
have been studied by several authors (\cite{CMc, H, HTU, Mc, MM}).
To state known results for type (ii) almost complex curves, we need
to recall Burstall and Pedit's definition of primitive maps
(\cite{BP}).  Let $\sigma$ be an order $6$ inner automorphism of
$\rmG_2$ such that the fixed point set of $\sigma$ is a maximal
torus $\rmT^2$, i.e., $\rmG_2/\rmT^2$ is a $6$-symmetric space. Let
$\fh_j$ denote the eigenspace of the complexified $d\sigma_e$ on
$\fg_2^\C=\fg_2\otimes \C$.  A map $f:\C\to \rmG_2/\rmT^2$ is {\it
primitive\/} if there is a lift $F:\C\to \rmG_2$ such that
$F^{-1}F_z\in \fh_0+\fh_{-1}$.   We will call any smooth map
$F:\C\to \rmG_2$ satisfying the condition that $F^{-1}F_z\in
\fh_0+\fh_{-1}$ a {\it $\sigma$-primitive $G_2$-frame\/}. Bolton,
Pedit, and Woodward (\cite{BPW}) proved that if $f:\Sigma\to S^6$ is
a type (ii) almost complex curve, then there exists a
$\sigma$-primitive $\rmG_2$-frame $\psi$.  Conversely, they show
that if $\psi$ is a $\sigma$-primitive $\rmG_2$-frame, then the
first column of $\psi$ gives an almost complex curve. The equation
for $\sigma$-primitive $\rmG_2$-frame is an elliptic integrable
system, so techniques from integrable systems can be used to study
almost complex surfaces in $S^6$.

In this paper, we prove that if $\Sigma$ is an immersed almost
complex surface in $S^6$ such that the second fundamental form $\II$
is not zero at $p_0$, then there exist an open neighbor $\co$ of
$p_0$ and a $\sigma$-primitive $\rmG_2$-frame $\psi:\co\to G_2$ such
that the first column is the immersion.  In other words, the
Gauss-Codazzi equation for the associative cones in $\R^7$ is the
equation for $\sigma$-primitive $\rmG_2$-frames.  Then we use this
elementary submanifold geometry set up to derive some of the known
properties of almost complex curves in $S^6$.  We also formulate the
equation for $S^1$-symmetric almost complex curves in $S^6$ as a
Toda type equation and use the AKS (Adler-Kostant-Symes) theory (cf.
\cite{AM, BP, AMV})  to construct $S^1$-symmetric almost complex
curves.

  This paper is organized as follows. We review basic properties of
$\rmG_2$ (\cite{Ko}) in section 2,  prove the existence of a
$\sigma$-primitive $\rmG_2$-frame on an almost complex surface with
non-vanishing second fundamental form in section 3. The equation for
$\sigma$-primitive $\rmG_2$-frame is a system of first order PDEs for
$5$ complex functions, we explain in section 4 the necessary and
sufficient conditions on these $5$ functions corresponding to the
four types of almost complex curves.
  In section 5, we explain how periodic Toda lattice arises from
$S^1$-symmetric almost complex curves in $S^6$, and finally in
section 6, we use the AKS theory to construct all $\rmS^1$-symmetric
almost complex curves.

\section{The octonions and Lie group $\rmG_2$}

Let $\mathbb{H} = \R \{{\bf 1}, {\rm {\bf i}}, {\rm {\bf j}} , {\rm
{\bf k}} \}$ be the quaternions, where ${\rm {\bf i}}$, ${\rm {\bf
j}}$ and  ${\rm {\bf k}}$ satisfy the condition ${\rm {\bf i}} \cdot
{\rm {\bf j}} = {\rm {\bf k}}, ~{\rm {\bf j}} \cdot {\rm {\bf k}} =
{\rm {\bf i}}, ~{\rm {\bf k}} \cdot {\rm {\bf i}} = {\rm {\bf j}},
~{\rm {\bf i}}^2 = {\rm {\bf j}}^2 = {\rm {\bf k}}^2 = -{\bf 1}$.
The conjugate of $a = a_0 + a_1 {\rm {\bf i}} + a_2 {\rm {\bf j}}
+a_3 {\rm {\bf k}}$ is $\overline{a} = a_0 - a_1 {\rm {\bf i}} - a_2
{\rm {\bf j}} - a_3 {\rm {\bf k}}$. The quaternions $\mathbb{H}$
equipped with the standard norm of $\R^4$ is an associative normed
algebra, i.e., $\parallel a \cdot b
\parallel =
\parallel a
\parallel \cdot \parallel b \parallel $. The octonions are defined
to be $ \O = \mathbb{H} \oplus \mathbb{H} {\rm {\bf e}} $ with the
multiplication
\[
( a + b {\rm {\bf e}} ) \cdot (c + d {\rm {\bf e}}) = ( a \cdot c -
\overline{d} \cdot b ) + ( d \cdot a + b \cdot \overline{c} ) {\rm
{\bf e}}
\]
The octonions $ \O$ equipped with the standard norm of $\R^8$ is a
non-associative normed algebra. Let $\{e_1, \cdots,e_7 \}$ be the
standard basis of $\R^7$. We identify $\R^7$ with ${\rm Im} \O$ as
follows:
\[
e_1 \rightarrow {\rm {\bf i}},~ e_2 \rightarrow {\rm {\bf j}}, ~ e_3
\rightarrow {\rm {\bf k}},~ e_4 \rightarrow {\rm {\bf e}},~ e_5
\rightarrow {\rm {\bf i}} {\rm {\bf e}},~ e_6 \rightarrow {\rm {\bf
j}} {\rm {\bf e}},~ e_7 \rightarrow {\rm {\bf k}} {\rm {\bf e}}.
\]
The multiplication table of octonions is:
\[
\begin{tabular}{|c|c|c|c|c|c|c|c|} \hline

     &  $e_1$  &  $e_2$   &  $e_3$  &  $e_4$  &  $e_5$  &  $e_6$  &
$e_7$       \CR \hline
  $e_1$   & $-1$  &  $e_3$   &  $-e_2$ & $e_5$  &  $-e_4$  & $-e_7$  &
$e_6$       \CR \hline
  $e_2$   & $-e_3$  &  $-1$  &  $e_1$  & $e_6$  & $e_7$   & $-e_4$   &
$-e_5$       \CR \hline
  $e_3$   &  $e_2$  &  $-e_1$  &  $-1$ & $e_7$  & $-e_6$  &  $e_5$  &
$-e_4$       \CR \hline
  $e_4$   & $-e_5$ & $-e_6$  & $-e_7$ & $-1$  &  $e_1$   &  $e_2$   &
$e_3$       \CR \hline
  $e_5$  &  $e_4$  &  $-e_7$ & $e_6$  & $-e_1$  &  $-1$  &  $-e_3$  &
$e_2$       \CR \hline
  $e_6$  & $e_7$  &  $e_4$   & $-e_5$ & $-e_2$  &  $e_3$   &  $-1$  &
$-e_1$       \CR \hline
$e_7$  & $-e_6$ &  $e_5$  &  $e_4$  & $-e_3$  &  $-e_2$  &  $e_1$ &
$-1$      \CR \hline
\end{tabular}
\]

The Lie group $\rmG_2$ is defined by
\[
\rmG_2 = {\rm Aut\/}(\O) = \{ g \in {\rm GL}(\O) ~|~ g(x \cdot y) =
g(x) \cdot g(y) \}
\]
We list below some basic properties of the Lie group $\rmG_2$ we
need in this paper: \ben
\item Let $f_1,f_2$ be two orthonormal column vectors in
$\R^7$. If $f_3 =f_1 \cdot f_2$, then $f_3$ is a unit vector and
perpendicular to $f_1,f_2$. Let $f_4$ be a unit column vector which
is perpendicular to $f_1,f_2,f_3$ and denote $f_5 = f_1 \cdot f_4,~
f_6 = f_2 \cdot f_4, ~f_7 = f_3 \cdot f_4$. Then $ (f_1, \cdots,f_7)
\in \rmG_2$ Such $\{f_1, \cdots, f_7 \}$ is called a $\rmG_2$-frame.
\item Any element of $\rmG_2$ can be realized by a $\rmG_2$-frame.
\item $\rmG_2$ is a compact, simply-connected, simple Lie group,
$\rmG_2 \subseteq \rmSO( {\rm Im} \O)$, and $\dim(\rmG_2) = 14$.
\item
Let $x^1, \cdots, x^7$ be coordinates of $\R^7$. The $3$-form
$\phi(x,y,z)= (x, y\cdot z)$ can be written as \[ \phi = d x^{123} +
d x^{145} - d x^{167} + d x^{246} - d x^{275} + d x^{347} - d
x^{356} \] where $d x^{jkl} = d x^j \wedge d x^k \wedge d x^l$. Then
\[ \rmG_2 = \{ ~g \in {\rm GL}(7, \mathbb{R}) ~|~ g^* \phi = \phi
~\} \]
  \item
The Lie algebra $\frak{g}_2$ of ${\rm G}_2$ are the space of
matrices \beq \label{g2} \bpm 0 & -x_2 & -x_3 & -x_4 & -x_5     &
-x_6 & -x_7 \cr x_2 & 0    & -y_3 & -y_4 & -y_5 & -y_6     & -y_7
\cr x_3 & y_3  & 0 & -x_6+y_5 & -x_7-y_4 & x_4-y_7  & x_5+y_6 \cr
x_4 & y_4  & x_6 - y_5 & 0 & -z_5     & -z_6     & -z_7 \cr x_5 &
y_5 & x_7 + y_4 & z_5 & 0 & -x_2-z_7 & -x_3+z_6 \cr x_6   & y_6 & -
x_4+y_7 & z_6 & x_2+z_7  & 0 &-y_3-z_5 \cr x_7   & y_7 & -x_5-y_6 &
z_7 & x_3-z_6 & y_3 +z_5 & 0   \epm \eeq where $x_2, \cdots, x_7$,
$y_3, \cdots,y_7$, $z_5,z_6,z_7$ are real numbers. To see this fact,
we let $\{ e_1, \cdots ,e_7 \}$ be the standard bases in $\R^7$. We
have $e_3 = e_1 \cdot e_2$, $e_5 = e_1 \cdot e_4$, $e_6 = e_2 \cdot
e_4$, $e_7 = ( e_1 \cdot e_2 ) \cdot e_4 $. If $ A \in \frak{g}_2$,
then
\[
A (e_j \cdot e_k ) = A ( e_j) \cdot e_k + e_j \cdot A (e_k)
\]
So $A$ is determined by $ A(e_1), A(e_2)$ and $A(e_4)$. Let $A(e_1)
= x_2 e_2 + \cdots + x_7 e_7$. Since $ A \in \frak{g}_2 \subset
\frak{so}(7) $, we can write $A(e_2) = - x_2 e_1 + y_3 e_3 + \cdots
+ y_7 e_7$. Then
\begin{align*} A(e_3) &= A ( e_1) \cdot e_2 + e_1 \cdot A (e_2) \\
       &= - x_3 e_1 - y_3 e_2 + ( x_6 -y_5) e_4 + ( x_7 + y_4 ) e_5\\
        &\quad + (y_7 -x_4 ) e_6 - ( x_5 + x_6) e_7
\end{align*}
  Since $A \in \frak{g}_2 \subset \frak{so}(7) $, we can write
\[
A (e_4) = - x_4 e_1 - y_4 e_2 + ( y_5 - x_6 )e_3 + z_5 e_5 +  z_6
e_6 + z_7 e_7
\]
Similarly $A(e_5),\cdots, A(e_7)$ are determined. Thus $A$ is a
matrix of type \eqref{g2}.  Conversely, any matrix of type
\eqref{g2} is a element of $\frak{g}_2$. \een

\section{$\sigma$-primitive $G_2$-frame}

Let $X_2$ denote the matrix defined by \eqref{g2} with $x_2 =1$, and
all other variables being zero. The matrices $X_3, \cdots, X_7$,
$Y_3, \cdots,Y_7$, $Z_5,Z_6,Z_7$ are defined similarly.

Let $h = {\rm exp}( \frac{\pi}{3} (Y_3 + 2 Z_5))$, and
$\sigma:\rmG_2\to \rmG_2$ the order $6$ inner automorphism defined
by $\sigma(g)= h^{-1}gh$.
  The eigenspace $\frak{h}_j$ with eigenvalue ${\rm exp}
\left({\frac{j\pi  {\rm i}}{3}} \right)$
  for the complexified $d\sigma_e$ on $\fg_2^{\C}=\fg_2\otimes \C$ is:
\beq \label{eigen}
\begin{array}{rcl}
\frak{h}_0 &=&\{ Y_3, Z_5 \} \CR \frak{h}_1 &=& \left\{ X_2 + {\rm
i} X_3 + \frac{{\rm i}}{2} \left(Z_6 + {\rm i} Z_7 \right), ~Y_4 +
{\rm i} Y_5, ~Z_6 - {\rm i} Z_7 \right\} \CR \frak{h}_2 &=& \left\{
X_4 + {\rm i} X_5 - \frac{{\rm i}}{2} \left(Y_6 + {\rm i} Y_7
\right), ~Y_6 - {\rm i} Y_7 \right\} \CR \frak{h}_3 &=& \left\{ X_6
- {\rm i} X_7 + \frac{{\rm i}}{2} \left(Y_4 - {\rm i} Y_5 \right),
~X_6 + {\rm i} X_7 - \frac{{\rm i}}{2} \left(Y_4 + {\rm i} Y_5
\right) \right\} \CR \frak{h}_4 &=& \left\{ X_4 - {\rm i} X_5 +
\frac{{\rm i}}{2} \left(Y_6 - {\rm i} Y_7 \right), ~Y_6 + {\rm i}
Y_7 \right\} \CR \frak{h}_5 &=& \left\{ X_2 - {\rm i} X_3 -
\frac{{\rm i}}{2} \left(Z_6 - {\rm i} Z_7 \right), ~Y_4 - {\rm i}
Y_5, ~Z_6 + {\rm i} Z_7 \right\}
\end{array}
\eeq Here $\{v_1, \ldots, v_m\}$ means the linear span of $v_1,
\ldots, v_m$. Notice $\bar\fh_j= \fh_{-j}$ (we use the convention
that $\fh_i=\fh_j$ if $i\equiv j$ (mod $6$)).

A smooth map $\psi:\C\to \rmG_2$ is $\sigma$-primitive if there exists $(u_0, u_{-1}):\C\to \frak{h}_0+\frak{h}_{-1}$ such that 
$$\psi^{-1} d\psi = (u_0+u_{-1}) dz + (\bar u_0 +\bar u_{-1}) d\bar z.$$
The flatness of $\psi^{-1}d\psi$ implies that $(u_0,u_{-1}):
\mathbb{C} \to \frak{h}_0 \oplus \frak{h}_{-1}$ must satisfy
\begin{equation}\label{aa}
\bca
(u_0)_{\bar{z}} - \left( \bar u_0 \right)_z = [u_0, \bar u_0] +
[u_{-1}, \bar u_{-1}],\\
(u_{-1})_{\bar{z}} = [u_{-1},\bar u_0]. \eca
\end{equation}
This system has a Lax pair \beq \label{Lax1} \theta_{\lambda} =
\left( u_0 +  {\lambda}^{-1} u_{-1} \right) \rmd z + \left( \bar u_0
+ \lambda \bar u_{-1} \right) \rmd \bar{z} \eeq i.e., $(u_0,u_{-1})$
is a solution of $(\ref{Lax1})$ if and only if $\theta_{\lambda}$ is
flat for all $\lambda \in \mathbb{C} \backslash \{0\}$. Note that:
\begin{itemize}
\item[(1)] The Lax pair  satisfies the  following reality conditions:
\beq\label{} \overline{ \left( \theta_{1/\bar{\lambda}} \right)} =
\theta_{\lambda}, \quad  \sigma(\theta_{\lambda}) =
\theta_{e^{\frac{\pi i}{3}} \lambda} \eeq
\item[(2)] $\xi(\lambda) = \sum_j \xi_j \lambda^j$ satisfies the
above reality condition if and only if $\xi_j \in \frak{h}_j$ and
$\xi_{-j} = \bar\xi_j$ for all $j$.
\end{itemize}

The following is well-known:

\bprop Let $(u_0, u_{-1}):\mathbb{C} \to  \frak{h}_0 \oplus
\frak{h}_{-1}$ be smooth maps. The following statements are
equivalent:
\begin{itemize}
\item[(1)] $(u_0, u_{-1})$ satisfies \eqref{aa}.
\item[(2)] $\theta_{\lambda} = \left( u_0 +  {\lambda}^{-1} u_{-1}
\right) \rmd z + \left( \overline{u}_0 + \lambda \overline{u}_1
\right) \rmd \bar{z}$ is flat for all $\lambda \in \mathbb{C}
\backslash \{0 \}$, i.e., $\rmd \o_\l= -\o_\l\wedge \o_\l$.

\item[(3)] $\theta_1 = \left( u_0 + u_{-1} \right) \rmd z + \left(
\overline{u}_0 +  \overline{u}_1 \right) \rmd \bar{z}$ is flat.

\item[(4)] There exists $\psi: \mathbb{C} \to \rmG_2$ such that
$\psi^{-1} \psi_z = u_0 + u_{-1}$, i.e., $\psi$ is a
$\sigma$-primitive $G_2$-frame.
\end{itemize}
\eprop

\begin{proof}
The only nontrivial part is $(3) \Leftrightarrow (1)$. To see this,
we decompose \beqn \rmd \theta + \theta \wedge \theta &=& \left(
-(u_0)_{\bar{z}} + \left( \overline{u}_0 \right)_z  + [u_{-1}, \bar
u_{-1}  \right) \rmd z \wedge \rmd \bar{z} \CR & & + \left( -
(u_{-1})_{\bar{z}} + [u_{-1},\overline{u}_0] \right) \rmd z \wedge
\rmd \bar{z} \CR & & + \left( \left(\bar  u_{-1} \right)_z +
[u_0,\bar u_{-1}] \right) \rmd z \wedge \rmd \bar{z} \eeqn according
to $\frak{h}_0 \oplus \frak{h}_1 \oplus \frak{h}_{-1}$. Thus
$(u_0,u_{-1})$ satisfies \eqref{aa} if and only if $\rmd \theta +
\theta \wedge \theta =0$.
\end{proof}

Suppose $(u_0, u_{-1})$ is a solution of \eqref{aa}. Since $\o_\l$
is flat at $\l=1$, there exists $\psi:\C\to G_2$ such that \beq
\label{matrix} \psi^{-1} \psi_z = u_0 + u_{-1} = \bpm 0 & - c & {\rm
i} c& & & & \cr c & 0 & - a & -d & {\rm i} d& & \cr
  -{\rm i} c & a & 0 & -  {\rm i} d & -d & & \cr
& d & {\rm i} d & 0 & -b & -e + \frac{{\rm i}}{2}c & - {\rm i} e +
\frac{1}{2} c \cr
  & -  {\rm i} d & d & b & 0 & - {\rm i} e -\frac{1}{2} c & e
+\frac{{\rm i}}{2} c \cr
  & & & e - \frac{{\rm i}}{2} c & {\rm i} e + \frac{1}{2} c & 0 &-a -b \cr
  & & & {\rm i} e -\frac{1}{2} c & -e - \frac{{\rm i}}{2} c & a+b & 0 \epm
\eeq System \eqref{aa} written in terms of $a, \ldots, e$ is \beq
\label{system} \bca
a_{\bar{z}} - \left( \overline{a} \right)_z = {\rm i} \left( 2 |c|^2
- 4 |d|^2 \right) \\
b_{\bar{z}} - \left( \overline{b} \right)_z = {\rm i} \left( - |c|^2
+ 4 |d|^2 - 4 |e|^2 \right) \\
c_{\bar{z}} = - {\rm i} \overline{a} c \\
d_{\bar{z}} = {\rm i} \left( \overline{a} - \overline{b} \right) d \\
e_{\bar{z}} = {\rm i} \left( \overline{a} + 2 \overline{b} \right) e
\eca \eeq Let $f_1, \ldots, f_7$ denote the columns of $\psi$. Then
\eqref{matrix} written in columns gives \beq\label{ao} \bca
(f_1)_z= cf_2-icf_3,\\
(f_2)_z= -c f_1 + af_3 + d\ f_4-id\  f_5,\\
(f_3)_z= ic f_1-a f_2 + id\ f_4+ d\ f_5,\\
(f_4)_z= -d\ f_2-id\ f_3 + bf_5 +(e-\frac{ic}{2}) f_6 +(ie-\frac{c}{2}) f_7,\\
(f_5)_z= id\ f_2 -d\ f_3 -bf_4 +(ie+\frac{c}{2})f_6-(e+\frac{ic}{2}) f_7,\\
(f_6)_z= (-e+\frac{i}{2}c) f_4 -(ie+\frac{c}{2}) f_5 +(a+b)f_7,\\
(f_7)_z= (-ie+\frac{c}{2})f_4 +(e+\frac{ic}{2})f_5 -(a+b)f_6. \eca
\eeq

\section{Associative cones and almost complex curves}

The following well-known  Proposition  relates almost complex curves
to associative cones:

\bprop (\cite{HL}) Let $\Sigma$ be a $2$-dimensional surface in
$\rmS^6$, and ${\rm C}(\Sigma) = \{tx ~|~ t > 0, x \in M \}$ the
cone of $\Sigma$ in $\R^7$. Then ${\rm C}(\Sigma)$ is an associative
submanifold in $\R^7$ if and only if $\Sigma$ is a almost complex
curve in $\rmS^6$. \eprop

\begin{proof}
Let $\{ e_1, e_2 \}$ be an orthonormal basis of $T_x \Sigma$. Then
$\{ x, e_1,e_2 \}$ is an orthonormal basis of $T_x {\rm C}(\Sigma)$.
Lemma follows from the fact that $\R \{ {\bf 1}, x, e_1,e_2 \}$ is
an associative subalgebra if and only if $x \cdot e_1 = e_2$.
\end{proof}

So the study of associative cones in $\R^7$  reduces to the study of
almost complex curves in $S^6$.

Since associative cones are calibrated by the $3$-form $\phi$, they
are minimal.  But a cone $C(\Sigma)$ in $\R^7$ is minimal if and
only if $\Sigma$ is minimal in $S^6$, so almost complex curves in
$S^6$ are minimal.

\bthm (\cite{BPW}) If $\psi = (f_1, \cdots, f_7) : \mathbb{C} \to
\rmG_2$ satisfies \beq \label{prim} \psi^{-1} \psi_z \in \frak{h}_0
\oplus \frak{h}_{-1} \eeq Then $f_1 : \mathbb{C} \to \rmS^6$ is
almost complex.   Conversely, if $f:\C \to S^6$ is a type (ii)
almost complex curve, i.e., $f$ is full and not totally isotropic,
then there exists a $\sigma$-primitive map $\psi:\C\to G_2$ such
that the first column of $\psi$ is $f$. \ethm

The first part of the above theorem is easy to see: Write
$\psi=(f_1, \ldots, f_7)$, and
\[
\psi^{-1} \psi_z = u_0 + u_{-1}.
\]
Then $u_0 + u_{-1}$ is given by $(\ref{matrix})$, so
\[
\left( f_1 \right)_z = c f_2 - {\rm i} c f_3
\]
By the definition of almost complex structure $J$ on $\rmS^6$, we
have
\[
J \left( f_1 \right)_z = f_1 \cdot \left( f_1 \right)_z = c f_3 +
{\rm i} c f_2 = {\rm i} \left( f_1 \right)_z
\]
So $f_1$ is almost complex.

\ms Next we prove that a $\sigma$-primitive $G_2$-frame exists on
any almost complex curve in $S^6$ with non-vanishing second
fundamental forms.

\bthm \label{al} Suppose $f_1:\Sigma\to S^6$ is an almost complex
curve such that the second fundamental form $\II$ is non-zero at
some $p_0 \in \Sigma$. Then there exists a neighborhood $\co$ of
$p_0$ and a $\sigma$-primitive $\rmG_2$-frame $\psi = \{f_1,
\cdots,f_7 \}$ on $\co$ such that  $f_2$ and $f_3$ are tangent to
the immersion, $\psi^{-1}\psi_z$  is given by \eqref{matrix} in terms of $5$ functions $a,
\ldots, e$, and  \eqref{system} is
the Gauss-Codazzi equation for $f_1$.  Moreover, the first and
second fundamental forms of $f_1$ are
\begin{align*}
&{\rm I}= 2|c|^2 |dz|^2,\\
&{\rm II} \left( \frac{\partial}{\partial z},
\frac{\partial}{\partial z} \right) = 2 cd (f_4 - {\rm i} f_5),
\end{align*}
and the normal connection is given by the lower $4\times 4$ matrices
\eqref{matrix}. \ethm

\begin{proof}
  Locally we
can choose orthonormal tangent frame $\{f_2,f_3\}$ such that $f_3 =
f_1 \cdot f_2$. Let $f_4$ be an arbitrary unit vector such that $f_4
~\bot~ {\rm span}_{\R} \{f_1,f_2,f_3\}$. Then we have a
$\rmG_2$-frame $\psi = \{f_1, \cdots, f_7 \}$ where $f_5 = f_1 \cdot
f_4, f_6 =f_2 \cdot f_4, f_7 = f_3 \cdot f_4$. Therefore we obtain a
$\frak{g}_2$-valued flat connection $1$-form $\omega = (\omega_{ij}
) = \psi^{-1} \rmd \psi$.

Write
\[
\rmd f_1 = f_2 \otimes \theta_2 + f_3 \otimes \theta_3,
\]
where $\theta_j$ is the dual $1$-form of $f_j$ for $j=2,3$.
Therefore
\[
\begin{array}{ll}
\omega_{21} = \theta_2, &\omega_{31} = \theta_3 \CR \omega_{\alpha
1} = 0, & 4 \leq \alpha \leq 7
\end{array}
\]
Since $\omega$ is $\frak{g}_2$-valued, we have
\[
\omega_{43} = - \omega_{52}, ~\omega_{53} = \omega_{42},
~\omega_{63} = \omega_{72}, ~\omega_{73} = - \omega_{62}
\]
Let
\[
\omega_{52} = a_2 \theta_2 + a_3 \theta_3, ~\omega_{62} = b_2
\theta_2 + b_3 \theta_3
\]
It follows from the  flatness of $(\w_{ij})$ that
\[
d \omega_{\alpha 1} + \sum_{j=1}^7 \omega_{\alpha j} \wedge
\omega_{j1} = 0,  ~(\alpha = 4,5),
\]
so we have \beqn (a_2 \theta_2 + a_3 \theta_3) \wedge \theta_2 +
\omega_{42} \wedge \theta_3 &=& 0 \CR \omega_{42} \wedge \theta_2 -
(a_2 \theta_2 + a_3 \theta_3) \wedge \theta_3 &=& 0 \eeqn Thus
\[
\omega_{53} = \omega_{42} = a_3 \theta_2 - a_2 \theta_3
\]
Similarly,
\[
\omega_{63} = \omega_{72} = b_3 \theta_2 - b_2 \theta_3
\]
Then the second fundamental form of immersion is given by \beqn {\rm
II} &=& \sum_{\alpha = 4}^7 f_{\alpha} \otimes \left( \omega_{\alpha
2} \otimes \theta_2 + \omega_{\alpha 3} \otimes \theta_3 \right) \CR
&=& v_1 \otimes (\theta_2 \otimes \theta_2 - \theta_3 \otimes
\theta_3) - v_2 \otimes (\theta_2 \otimes \theta_3 + \theta_3
\otimes \theta_2 ) \eeqn where $v_1= a_3 f_4 + a_2 f_5 + b_3 f_7 +
b_2 f_6$ and $v_2 = a_2 f_4 - a_3 f_5 + b_2 f_7 - b_3 f_6$. Note
that
\[
(v_1, v_1)= (v_2, v_2), \quad (v_1, v_2) =0
\]
Since ${\rm II}(p_0) \neq 0$, there exists a neighborhood $U$ of $p$
such that $v_1$ and $v_2$ are nonzero. Let $\widetilde{f}_j = f_j, ~j
=1,2,3,$
\[
\widetilde{f}_4 = \frac{v_1}{||v_1||}
\]
and
\[
\widetilde{f}_5 = \widetilde{f}_1 \cdot \widetilde{f}_4, \quad
\widetilde{f}_6 = \widetilde{f}_2 \cdot \widetilde{f}_4, \quad
\widetilde{f}_7 = \widetilde{f}_3 \cdot \widetilde{f}_4
\]
Then $\widetilde{\psi} = \{ \widetilde{f}_1, \cdots, \widetilde{f}_7
\}$ is a $\rmG_2$-frame, and a computation using the octonion
multiplication implies that $\ti f_5= v_2/||v_2||$.  Let $\widetilde{\omega}
= \left( \widetilde{\omega}_{ij} \right) = \widetilde{\psi}^{-1}
\rmd \widetilde{\psi}$. Since $(\II, \ti f_6)=(\II, \ti f_7)=0$, we
have
\[
\widetilde{\omega}_{62} = \widetilde{\omega}_{63} =
\widetilde{\omega}_{72} = \widetilde{\omega}_{73} \equiv 0.
\]
So $\ti \omega$ lies in $\fh_0+\fh_1+\fh_{-1}$, where $\fh_j$ is the
eigenspace of $d\sigma$ on $\fg_2\otimes \C$ with eigenvalue
$e^{\frac{2\pi ji}{6}}$.  Or equivalently, $\psi^{-1}\psi_z$ is of
the form \eqref{matrix}, i.e., $\psi$ is a $\sigma$-primitive
$G_2$-frame.   In particular, this shows that the Gauss-Codazzi
equation for almost complex curves is \eqref{system}.  It follows
from \eqref{ao} and a computation that the two fundamental forms
for $f_1$ are given as in the Theorem.
\end{proof}

As a consequence of the Fundamental Theorem of submanifolds in space
forms and the above theorem, we get

\bcor \label{Curve} Every simply connected immersed almost complex
curve in $\left( \rmS^6, J \right)$ with non-vanishing second
fundamental form has a $\sigma$-primitive $G_2$-frame such that the
first column is the immersion. Conversely, the  first column of  a
$\sigma$-primitive $G_2$-frame  is an almost complex  surface in
$\rmS^6$. \ecor

Next, we use Theorem \ref{al} to give conditions on $a, \ldots, e$
to determine the four types of almost complex curves mentioned in
the introduction.

\bcor \label{ap} Let $(a, \ldots, e)$ be a solution of
\eqref{system}, $\psi$ a solution of \eqref{matrix}, and $f_1$ the
first column of $\psi$.  Then  $f_1$  is almost complex in $S^6$ and
is
\begin{itemize}
\item[(i)] full in $S^6$ and totally isotropic  if and
only if $e \equiv 0$ and $d \neq 0$,
\item[(ii)] full in $S^6$ and not totally isotropic if and only if $de\neq 0$,
\item[(iii)]  full in $\rmS^5$ if and only if
$de\neq 0$ and $a+b \equiv 0$,
\item[(iv)] totally geodesic two sphere if and only if $d \equiv 0$,
i.e., ${\rm II}\equiv 0$.
\end{itemize}
Moreover, the cone over the curve of type (iii) is a special
Lagrangian cone in $\R^6$ with the appropriate complex structure.
\ecor

\begin{proof}
The first fundamental form is positive definite, so $c\not=0$.  A
surface is full then $\II$ can not be zero, so $d\not=0$. Let $\psi$
satisfy $\psi^{-1}d\psi= (u_0+u_{-1})dz+ (\bar u_0+\bar u_{-1})d\bar
z$, and $f_1$ denote the first column of $\psi$, where
$u_0+u_{-1}\in \fh_0+\fh_{-1}$ is given by \eqref{matrix}.  Then
$f_1$ is almost complex. Use \eqref{ao} and a direct computation to
see that
\[
  ((\K_{\frac{\p}{\p z}})^2f_\ast(\frac{\p}{\p z}),\ \
(\K_{\frac{\p}{\p z}})^2f_\ast(\frac{\p}{\p z})) = -32 {\rm i} c^3
d^2 e,
\]
where $(Y,Z)= \sum_j y_j z_j$ is the complex bilinear form on $\C^7$.
If $f$ is totally isotropic, then
$$  ((\K_{\frac{\p}{\p z}})^if_\ast(\frac{\p}{\p z}),\ \
(\K_{\frac{\p}{\p z}})^jf_\ast(\frac{\p}{\p z}))=0$$
for all other $ 0 \leq i,  j \leq 2$, so  $e = 0$.

Next we prove that if an almost complex curve is of type (iii), then
$a+b \equiv 0$. Since there is a constant unit normal vector field
on the curve, there exists real functions $\lambda_i$ ($4 \leq i
\leq 7$) on the curve such that this normal vector is $\sum_{i=4}^7
\lambda_i f_i$. Then
\begin{align*}
  & ( \sum_{i=4}^7 \lambda_i f_i )_z = \sum_{i=4}^7 (\lambda_i)_z f_i +
\sum_{i=4}^7 \lambda_i (f_i)_z  \\
    &= \sum_{i=4}^7 (\lambda_i)_z f_i + \lambda_4 [-d f_2 - id f_3 + b
f_5 + (e-ic/2)f_6 + (ie-c/2)f_7 ]  + \cdots =0.
\end{align*}
So the coefficient of $f_i$  must be zero for $2\leq i\leq 7$. Since
$d \neq 0$, it implies that $\lambda_4 =0$ and $\lambda_5 =0$. The
coefficients for $f_6$ and $f_7$ are $(\lambda_6)_z - (a+b) $ and
$(\lambda_7)_z + (a+b)$ respectively. Therefore
$(\lambda_6+\lambda_7)_z=0$, i.e., $\lambda_6+\lambda_7$ is
anti-holomorphic. Since $\lambda_6+\lambda_7$ is also real, it must
be a constant. Finally both $\lambda_6$ and $\lambda_7$ have to be
constant because their square sum is $1$. Thus
$a+b=(\lambda_6)_z=0$.

Conversely, if $a+b \equiv 0$, then the system $(\ref{system})$
implies that
\[
c_{\bar{z}} = - {\rm i} \overline{a} c, \quad e_{\bar{z}} = - {\rm
i} \overline{a} e
\]
and ${\rm i} \left( |c|^2 - 4 |e|^2 \right)  = (a+b)_{\bar{z}} -
\left( \overline{a} + \overline{b} \right)_z = 0$. Let $\alpha =
\frac{c}{2e}$. Then $\alpha_{\bar{z}} =0$ and $|\alpha|=1$. So
$\alpha \in \rmS^1$ is a constant and $\beta = \frac{-1 + {\rm i}
\alpha}{- {\rm i} + \alpha}$ is a real constant. It follows from
\eqref{ao} that $\left( f_6 - \beta f_7 \right)_{\bar{z}} = 0$. Thus
$ n = \frac{1}{\sqrt{1 + \beta^2}} \left( f_6 - \beta f_7 \right)$
is a unit constant normal vector. So the image of the immersion lies
in the hyperplane $V$ which is orthogonal to $n$. Note $J(x) = n
\cdot x$ defines a complex structure on the hyperplane $V$ and $J
(f_1) = \frac{1}{\sqrt{1 + \beta^2}} ( \beta f_6 + f_7)$, $J (f_2) =
\frac{1}{\sqrt{1 + \beta^2}} ( f_4 + \beta f_5)$, and $J (f_3) =
\frac{1}{\sqrt{1 + \beta^2}} ( - \beta f_4 - f_5)$.  Thus $J
(\rmspan_{\R} \{f_1, f_2,f_3 \}) = \rmspan_{\R} \{ f_4, f_5, \beta
f_6 + f_7 \}$, so the cone over the image of $f_1$ is Lagrangian in
$\left( \R^6, J \right)$. We know it is minimal, so by Proposition
$2.17$ of \cite{HL} that it is $\theta$-special Lagrangian for some
$\theta$.
\end{proof}

Next we  use Theorem \ref{al} to give a proof of one of Bryant's
results on almost complex curves in $\rmS^6$.  First recall that the
$5$-dimensional complex quadric $Q_5$ is defined by
\[
Q_5 = \{[z_1 : \cdots : z_7 ] \in \mathbb{C} {\rm P}^6 ~|~ z_1^2 +
\cdots + z_7^2 = 0 \}.
\]

\bthm \cite{B0} If $f: \Sigma \to \rmS^6$ is a totally isotropic
almost complex curve that is not totally geodesic, then it can be
lifted to a horizontal holomorphic map to $Q_5$. \ethm

\begin{proof} Let  $\psi = (f_1, \cdots,
f_7) : \Sigma \to {\rm G}_2$ denote the $\sigma$-primitive
$G_2$-frame obtained in Theorem \ref{al}. So $\psi^{-1}\psi_z$ is of
the form \eqref{matrix}.  Let $\Phi : \Sigma \to Q_5 $ be the map
defined by
\[
\Phi = [f_6 + i f_7]
\]
Clearly $\Phi$ is well-defined and is independent of choice of the
frame. By \eqref{ao}, we have
\[
(f_6 + {\rm i} f_7)_{\bar{z}} = -2 \bar{e} (f_4 - {\rm i} f_5) -
{\rm i} (\bar{a} + \bar{b}) (f_6 + {\rm i} f_7 )
\]
But we have shown in Corollary \ref{ap} that if $f$ is totally
isotropic then $e=0$, so  $\Phi$ is holomorphic.
\end{proof}

\section{$\rmS^1$-symmetric solutions and periodic Toda lattice}

By the maximal torus theorem, given $A\in \cg_2$, there exists $k\in
\rmG_2$ and real numbers $\l_1, \l_2$  such that
$A=k^{-1}(\l_1Y_3+\l_2Z_5)k$. Note \beqn \l_1 Y_3+\l_2 Z_5= \bpm 0 &
& & & & & \cr & & -\lambda_1 & & & & \cr & \lambda_1 & & & & & \cr &
& & & -\lambda_2 & & \cr & & & \lambda_2 & & & \cr & & & & &
&\lambda_3 \cr & & & & & -\lambda_3 & \epm \eeqn where
$\l_3=-(\l_1+\l_2)$. We say $A=k^{-1}(\l_1Y_3+\l_2Z_5)k$ is {\it
rational\/} if $\l_1,\l_2$ are linearly dependent over the
rationals.  It is easy to see that $A$ is rational if and only if
$\{\exp(sA)\n s\in \R\}$ is periodic.

To construct a $S^1$-symmetric almost complex curve in $\rmS^6$, we
need to construct $\psi = e^{A s} g(t)$ with rational $A$ and $g(t)
\in \rmG_2$ such that
$$\psi^{-1}\psi_z= u_0+u_{-1} \in \frak{h}_0+\frak{h}_1,$$
  where $z=s+it$ and $u_0 + u_{-1}$ is given by $(\ref{matrix})$ and
$a,b,c,d, e$ are complex valued functions of $t$ only.   A simple
computation gives
\[
\psi^{-1} \rmd \psi = \left( g^{-1} Ag \right) \rmd s + \left(
g^{-1} g_t \right) \rmd t.
\]
The flatness of $\psi^{-1} \rmd \psi$ implies that \beq \label{ab}
\left( g^{-1} A g \right)_t = \left[g^{-1} Ag , g^{-1}g_t \right].
\eeq Write $a = a_1 + {\rm i} a_2, ~b = b_1 + {\rm i} b_2, \ \
\ldots, \ \  e= e_1+{\rm i} e_2$ in real and imaginary part, and $c
= r_1 e^{{\rm i} \beta_1}, ~d = r_2 e^{{\rm i} \beta_2}, ~e = r_3
e^{{\rm i} \beta_3}$ in polar coordinates. Since $\psi^{-1}\psi_s=
g^{-1}A g= \psi^{-1}\psi_z + \psi^{-1}\psi_{\bar z}$,
$\psi^{-1}\psi_t= g^{-1}g_t= {\rm i}\ (\psi^{-1}\psi_z -
\psi^{-1}\psi_{\bar z})$, and $\psi^{-1}\psi_{z}$ is given by
\eqref{matrix}, we have
\begin{align*}
  &g^{-1} A g =\bpm 0 & - 2c_1 & -2 c_2& & & & \cr
2 c_1 & 0 & - 2 a_1 & -2 d_1 & -2 d_2 & & \cr
  2 c_2 & 2 a_1 & 0 & 2 d_2 & -2 d_1 & & \cr
& 2 d_1 & -2 d_2 & 0 & - 2b_1 & -2 e_1 -c_2 & 2 e_2 + c_1 \cr
  & 2 d_2 & 2 d_1 & 2 b_1 & 0 & 2 e_2 -c_1 & 2 e_1 - c_2 \cr
  & & & 2 e_1 + c_2 & -2 e_2 + c_1 & 0 &-2 a_1 -2 b_1 \cr
  & & & -2 e_2 -c_1 & -2 e_1 + c_2 & 2 a_1+ 2 b_1 & 0\epm,\\
  &g^{-1} g_t=\bpm 0 & 2c_2 & -2 c_1& & & & \cr
-2 c_2 & 0 &  2 a_2 & 2 d_2 & -2 d_1 & & \cr
  2 c_1 & -2 a_2 & 0 & 2 d_1 & 2 d_2 & & \cr
& -2 d_2 & -2 d_1 & 0 & 2b_2 & 2 e_2 -c_1 & 2 e_1 - c_2 \cr
  & 2 d_1 & - 2 d_2 & - 2 b_2 & 0 & 2 e_1 + c_2 & -2 e_2 - c_1 \cr
  & & & -2 e_2 + c_1 & -2 e_1 - c_2 & 0 & 2 a_2 + 2 b_2 \cr
  & & & -2 e_1 + c_2 & 2 e_2 + c_1 & - 2 a_2- 2 b_2 & 0 \epm.
\end{align*}
System \eqref{ab} written in $a, b, r_i, \b_i$ gives the following
two separable systems
$$ \bca
\dot{a}_1 = 2 r_1^2 - 4 r_2^2,\\
\dot{b}_1 = - r_1^2 + 4 r_2^2 - 4 r_3^2,\\
\dot{r}_1 = - 2 a_1r_1,\\
\dot{r}_2 = 2 (a_1 -b_1)r_2,\\
\dot{r}_3= 2 (a_1 + 2b_1) r_3, \eca\qquad
 \bca
\dot{\beta}_1 = 2 a_2,\\
\dot{\beta}_2 =- 2 a_2 + 2 b_2,\\
\dot{\beta}_3 = -2 a_2 - 4b_2. \eca$$ So we may assume that
$a_2=b_2=\b_1=\b_2=\b_3=0$, i.e.,
$$a_2=b_2=c_2=d_2=e_2=0.$$
Substitute these conditions to the matrix formulas for $g^{-1}Ag$
and $g^{-1}g_t$ to  get
\begin{align*}
&P:=g^{-1} A g =\bpm 0 & - 2c_1 & & & & & \cr 2 c_1 & 0 & - 2 a_1 &
-2 d_1 &  & & \cr
   & 2 a_1 & 0 &  & -2 d_1 & & \cr
& 2 d_1 &  & 0 & - 2b_1 & -2 e_1  &  c_1 \cr
  &  & 2 d_1 & 2 b_1 & 0 &  -c_1 & 2 e_1  \cr
  & & & 2 e_1 &  c_1 & 0 &-2 a_1 -2 b_1 \cr
  & & & -c_1 & -2 e_1  & 2 a_1+ 2 b_1 & 0\epm,\\
&Q:= g^{-1} g_t=\bpm 0 &  & -2 c_1& & & & \cr & 0 &   & & -2 d_1 & &
\cr
  2 c_1 &  & 0 & 2 d_1 &  & & \cr
&  & -2 d_1 & 0 &  &-c_1 & 2 e_1  \cr
  & 2 d_1 &  &  & 0 & 2 e_1 & - c_1 \cr
  & & &  c_1 & -2 e_1 & 0 &  \cr
  & & & -2 e_1  &  c_1 &  & 0 \epm
\end{align*}
Since $\psi^{-1}\psi_z= u_0+u_{-1}\in \frak{h}_0+\frak{h}_{-1}$, $P=
u_0+\bar u_0+ u_{-1}+\bar u_{-1}$ and $Q=-i(u_0-\bar u_0+
u_{-1}-\bar u_{-1})$. By assumption that $a, b, \ldots, e$ are real,
so $u_0=\bar u_0$, and \beq\label{aw} P= 2u_0+ u_{-1} +\bar u_{-1},
\quad Q= i(u_{-1}-\bar u_{-1}), \eeq where
$$\bca
u_0 =a_1Y_3+b_1 Z_5\ \ \in \frak{h}_0\cap \fg_2, \\
u_{-1} = c_1(X_2-\frac{Z_7}{2}) + i(X_3+\frac{Z_6}{2}) + d_1 (Y_4+
iY_5) + e_1 (Z_6 - iZ_7)\ \  \in \frak{h}_{-1}. \eca$$
   Thus we have

\bprop Suppose $(u_0, u_{-1}):\R\to (\fh_0\cap \fg_2)\times
\fh_{-1}$ satisfies \beq \label{at} (2u_0+ u_{-1} + \bar u_{-1})_t=
[2u_0+ u_{-1} + \bar u_{-1}, \ i(u_{-1}-\bar u_{-1})], \eeq and
there exist a constant $A\in (\fh_0\cap\fg_2)+ \fh_{-1}$ and
$g:\R\to \rmG_2$ such that \beq\label{au} \bca
g^{-1}Ag= 2u_0+ u_{-1}+\bar u_{-1}, \\
g^{-1}g_t= u_{-1} -\bar u_{-1}. \eca\eeq Then $f(s,t)= e^{As} g(t)$
is an almost complex curve in $S^6$.  Moreover, $f$ is
$S^1$-symmetric if and only if $A$ is rational, and is doubly
periodic if and only if $A$ is rational and $g$ is periodic. \eprop

Define $v_1,v_2,v_3$ by
$$\bca e^{2 v_1}=c_1^2\\ e^{2 (v_2 - v_1)}=d_1^2 \\
  e^{2 (v_3 -v_2)}= e_1^2 . \eca$$
  Then $a_1, b_1, v_1, v_2, v_3$ satisfy
\beq \label{aq} \bca
\dot{a}_1= 2 e^{2 v_1} - 4 e^{2 (v_2 - v_1)},\\
\dot{b}_1 = - e^{2 v_1} + 4 e^{2 (v_2 - v_1)} - 4 e^{2 (v_3 -v_2)},\\
\dot{v}_1= -2 a_1,\\
\dot{v}_2 = -2 b_1,\\
\dot{v}_3 = 2 (a_1 +b_1). \eca \eeq Clearly, $(v_1+v_2+v_3)_t =0$.
Moreover, $v_1,v_2,v_3$ satisfy
$$\bca
\ddot{v}_1 = -4 e^{2 v_1} + 8 e^{2 (v_2 - v_1)},\\
\ddot{v}_2 = 2 e^{2 v_1} - 8 e^{2 (v_2 - v_1)} + 8 e^{2 (v_3 -v_2)},\\
\ddot{v}_3 = 2 e^{2 v_1} - 8 e^{2 (v_3 -v_2)}. \eca$$ These are
equivalent to the periodic Toda lattice equations of $\rmG_2$-type.

  If $a_1 + b_1 =0$, i.e., the type (iii) case, then $\dot a_1+\dot
b_1= e^{2v_1} - 4e^{2(v_3-v-2)}=0$, $\dot v_1+\dot v_2=\dot v_3=0$,
so there is  a positive constant $C_1$ such
that:
\[
e^{2 (v_1 + v_2)} = 4 e^{2 v_3} = C_1 .
\]
  Then $v_1$ satisfies
\[
\ddot{v}_1 + 4 e^{2v_1} - 8 C_1 e^{-4 v_1} = 0
\]
Multiply $\dot{v}_1$ to both sides and integrating once to get
\[
\left( \dot{v}_1 \right)^2 + 4 e^{2v_1} + 4 C_1 e^{-4 v_1} = 4 C_2,
\]
where $C_2$ is a positive constant. Let $y = e^{2 v_1} = r_1^2$.
Then the above equation becomes
\[
\left( \dot{y} \right)^2 = -16 y^3 + 16 C_2 y^2 -16 C_1.
\]
One can verify easily that $4 C_2^3 \geq 27 C_1$. Therefore this
equation has three real constant solutions $\Gamma_1, \Gamma_2,
\Gamma_3$. Let us label these solutions so that $\Gamma_1 < 0 <
\Gamma_2 \leq \Gamma_3$. Then we can rewrite the previous equation
as
\[
\left( \dot{y} \right)^2 = -16 \left( y - \Gamma_1 \right)\left( y -
\Gamma_2 \right) \left( y - \Gamma_3 \right)
\]
Haskins (\cite{H}) showed that this  equation has the following
solution:
\[
y = \Gamma_3 - \left( \Gamma_3 - \Gamma_2 \right) {\rm sn}^2( B_1 t
+ B_2, B_3 )
\]
where $B_2$ is a constant determined by the initial condition of
$y$,
\[
B_1^2 = 4 (\Gamma_3 - \Gamma_1), \qquad B_3^2 = \frac{\Gamma_3 -
\Gamma_2}{\Gamma_3 - \Gamma_1}
\]
and sn is the Jacobi elliptic sn-noidal function. Recall that ${\rm
sn}(t, k)$ is defined to be the unique solution of the equation
\[
\dot{z}^2 = \left( 1 -z^2 \right)\left( 1 - k^2 z^2 \right)
\]
with $z(0) =0, \dot{z}(0) =1$, where $ 0 \leq k \leq 1$. It is
straightforward to see from this definition that ${\rm sn}(t,0) =
{\rm sin} ~t$ and ${\rm sn}(t,1) = {\rm tanh} ~t$. The period of
${\rm sn}(t, k)$ is given by
\[
\int_0^{2 \pi} \frac{dx}{\sqrt{1 - k^2 {\rm sin}^2 x}}
\]
Thus $y$ is a periodic function, so are $a_1,b_1,v_1,v_2$. They all
have same period denoted by ${\rm T}$.

In fact,  Haskins proved in \cite{H} that not only \eqref{at} has a
periodic solution but he also proved that the solution $g$ of
\eqref{au} is also periodic for some rational $A$.  So he proved the
existence of infinitely many $S^1$-symmetric type (iii) almost
complex curves (hence infinitely many special Lagrangian cones in
$\C^3$).

\section{$S^1$-symmetric solutions and loop group factorization}

The first equation of \eqref{au} implies that the solution $2u_0(t)+
u_{-1}(t)+ \bar u_{-1}(t)$ must lie in the same conjugate class for all
$t$, and there is $g$ solves \eqref{au}.  Although these conditions
seem to be extra conditions for solutions of \eqref{at}, we will see
below that \eqref{at} has a Lax pair and is a Toda type equation,
and hence the AKS theory implies that if $(u_0, u_{-1})$ is a
solution of \eqref{at} then there exists $g$ satisfies \eqref{au}
automatically.

Set $P= 2u_0 + u_{-1} + \bar u_{-1}$ and $Q= i(u_{-1}- \bar u_{-1})$
as in \eqref{aw}.  Then \eqref{au} is $P_t=[P,Q]$, or equivalently,
$iP_t= [iP, Q]$, i.e., \beq\label{af} (v_0+ v_{-1} -\bar
v_{-1})_t=[v_0+ v_{-1} -\bar v_{-1}, v_{-1}+\bar v_{-1}], \eeq where
$v_0\in \frak{h}_0\cap i\fg_2$ and $v_{-1}\in \frak{h}_{-1}$.

\ms \ni {\bf Equation \eqref{af} has a Lax pair}

\ss
   A simple calculation shows that $(v_0, v_1)$ satisfies \eqref{af}
if and only if
   \beq\label{ah}
   (v_0+v_{-1}\l^{-1}-\bar v_{-1}\l)_t= [v_0+v_{-1}\l^{-1}-\bar
v_{-1}\l, \ v_{-1}\l^{-1}+ \bar v_{-1}\l]
   \eeq
   holds for all $\l\in \C\setminus \{0\}$.  Here $v_0\in \frak{h}_0$
is pure imaginary and $v_{-1}\in \frak{h}_{-1}$.

\ms \ni {\bf  Results from the Adler-Kostant-Symes (AKS) Theory (cf.
\cite{AM, BP, AMV})}

\ss Let $G$ be a group, $G_+, G_-$ subgroups of $G$ such that the
multiplication map $G_+\times G_-\to G$ defined by $(g_+, g_-)\to
g_+g_-$ is a bijection.  So $\cg= \cg_++\cg_-$ as direct sum of
vector subspaces.  Suppose $\cg$ admits a non-degenerate,
ad-invariant bilinear form $(\ , )$.  Let \beq\label{aj}
\cg_+^\perp=\{y\in \cg\n (y, x)=0\ \forall \ x\in \cg_+\}, \eeq
  and $\pi_+$ denote the projection of $\cg$ onto $\cg_+$ with respect
to the decomposition $\cg=\cg_++\cg_-$.
  Suppose $M\subset \cg_+^\perp$ is invariant under the flow
$$\frac{dx}{dt}= [x(t), \pi_+(x(t))].$$
  Given $x_0\in M$, consider the following ODE:
\beq\label{ac}
\bca \frac{dx}{dt}= [x(t), \pi_+(x(t))],\\
x(0)=x_0. \eca \eeq The AKS theory gives a method to solve the
initial value problem \eqref{ac} via factorizations as follows:
\begin{itemize}
\item[(i)] Find the one-parameter subgroup $f(t)$ for $x_0$, i.e.,
solve $f^{-1}f_t= \ x_0$ with $f(0)=e$.
\item[(ii)] Factor $f(t)= f_+(t) f_-(t)$ with $f_\pm(t)\in G_\pm$.
\item[(iii)] Set $x(t)= f_+(t)^{-1} x_0 f_+(t)$.  Then $x(t)$ is the
solution for the initial value problem \eqref{ac}.  Moreover,
$f_+^{-1}(f_+)_t= \pi_+(x(t))$.
\end{itemize}
If $G=SL(n,\R)$, $G_+=SO(n)$, $G_-=$ the subgroup of upper
triangular matrices, and $M$ is the space of all tri-diagonal
matrices in $sl(n,\R)$, then ODE \eqref{ac} is the standard Toda
lattice.  So we call a system obtained from a factorization a {\it
Toda type equation\/}.

   \ms
   \ni {\bf Equation \eqref{af} is of Toda type}

   Let $L(\rmG_2^\C)$ denote the group of smooth loops from $S^1$ to
$\rmG_2^\C$ satisfying the reality condition
$\overline{g(\bar\l^{-1})}= g(\l)$,
$L_+(\rmG_2^\C)$ the subgroup of $g\in L(\rmG_2^\C)$  with $g(\l)\in
\rm G_2$ for all $\l\in S^1$, and $L_-(\rmG_2^\C)$ denote the
subgroups of $f\in L(\rmG_2^\C)$ that can be extended to a
holomorphic maps inside $S^1$ such that $f(0)=e$ the identity of
$G$.    Pressely and Segal proved in \cite{PS} an analogue of the
Iwasawa decomposition of simple Lie groups for  loop groups:

\bthm {\rm (Iwasawa loop group factorization Theorem \cite{PS,
G})}\par The multiplication map $L_+(\rmG_2^\C)\times
L_-(\rmG_2^\C)\to L(\rmG_2^\C)$ is a diffeomorphism. In particular,
given $g\in L(\rmG_2^\C)$, we can factor $g=g_+g_-$ uniquely with
$g_\pm \in L_\pm(\rmG_2^\C)$.
  \ethm

  Note that
$$\hat \sigma(g)(\l)= \sigma(g(e^{-\frac{\pi i}{3}}\l))$$
defines an automorphism of $L(\rmG_2^\C)$.  Let
$L^\sigma(\rmG_2^\C)$ and $L_\pm^\sigma(\rmG_2^\C)$ denote the
subgroups fixed by $\hat \sigma$ of $L(\rmG_2^\C)$ and
$L_\pm(\rmG_2^\C)$ respectively.   Then we have

\bcor If  $g\in L^\sigma(\rmG_2^\C)$ and $g=g_+g_-$ with $g_\pm \in
L_\pm(\rmG_2^\C)$, then $g_\pm \in L^\sigma_\pm(\rmG_2^\C)$. \ecor

Let $B$ denote the Borel subgroup of $\rmG_2^\C$ such that the
Iwasawa decomposition is $\rmG_2^\C= G_2B$, and $\fg_2^\C=
\fg_2+\fb$ at the Lie algebra level. It is easier to write down the
factorization at the Lie algebra level: \beq\label{av}
\cl^\sigma(\fg_2^\C)= \cl^\sigma_+(\fg_2^\C)
+\cl_-^\sigma(\fg_2^\C), \eeq where
\begin{align*}
&\cl^\sigma(\fg_2^\C)= \{\xi=\sum_{j\in {\Z}} \xi_j\l^j\n \xi_j\in
\fg_2^\C, \xi_j\in \fh_j\},\\
&\cl^\sigma_+(\fg_2^\C)=\{\xi=\sum_{j\in {\Z}} \xi_j\l^j\in
\cl^\sigma(\fg_2^\C)\n \xi_{-j}= \bar \xi_j\},\\
&\cl^\sigma_-(\fg_2^\C)= \{\xi=\sum_{j\geq 0}
\xi_j\l^j\in\cl^\sigma(\fg_2^\C)\n \xi_0\in \fb\}.
\end{align*}
Let $\pi_{\fg_2}$ and $\pi_{\fb}$ denote the projections of
$\fg_2^\C$ onto $\fg_2$ and $\fb$ respectively, and $\pi_\pm$ the
projections of $\cl^\sigma(\fg_2^\C)$ onto
$\cl^\sigma_\pm(\fg_2^\C)$ with respect to the decomposition
\eqref{av}.  Then for $\xi= \sum_j \xi_j \l^j$,
\begin{align*}
\pi_+(\xi)&= \pi_{\fg_2}(\xi_0)+ \sum_{j> 0}
\xi_{-j}\l^{-j}+\bar\xi_{-j}\l^j, \\
\pi_-(\xi)&= \pi_{\fb}(\xi_0) +\sum_{j> 0} (\xi_j-\bar\xi_{-j})\l^j.
\end{align*}

  Let $(\ , )$ be the Killing form on $\cg_2^{\C}$. Then
   $$\li \xi, \eta\ri =\sum_{i+j=0} (\xi_i, \eta_j)$$
   is an ad-invariant bilinear form on $\cl(\cg)$.  So
   $$\cl_+(\cg)^\perp=\{\xi=\sum_j\xi_j\l^j\n \xi_{-j}= -\bar\xi_j\}.$$
Let $M=\{\xi= \xi_0 + \xi_{-1}\l^{-1} -\bar\xi_{-1}\l \n \xi_0\in
\frak{h}_0\cap (i\cg_2), \xi_{-1}\in \frak{h}_{-1}\}$.  Note that
   $$\pi_+(\xi_0+\xi_{-1}\l^{-1} - \bar\xi_{-1}\l)= \xi_{-1}\l +
\bar\xi_{-1} \l.$$
   It is easy to check that $[\xi, \pi_+(\xi)]\in M$ if $\xi\in M$, so
$M$ is invariant under the flow $\xi_t=[\xi, \pi_+(\xi)]$.   So we
can use the Iwasawa loop group factorization to construct solution of
\eqref{ah} as described in the AKS theory and get

   \bthm
   Let $A= 2h_0+h_{-1}+ \bar h_{-1}$ be a constant with $h_0\in
\fh_0\cap \fg_2$ and $h_{-1}\in \fh_{-1}$.  Then the solution of
\eqref{at} with initial value $A$ can be obtained as follows:
   \ben
   \item  Set $\xi_0(\l)= 2ih_0 + ih_{-1}\l^{-1} +i\bar h_{-1}\l$, and
construct $g(t,\l)$ such that
   $$\bca g^{-1} g_t= \xi_0(\l),\\ g(0,\l)= \I,
   \eca$$
i.e., $g(t, \cdot)$ is the one-parameter subgroup of $\xi_0$ in
$L^\sigma(\rmG_2^\C)$.
\item Factor $g(t,\l)= g_+(t, \l) g_-(t, \l)$ such that $g_\pm(t,
\cdot)\in L^\sigma_\pm(\rmG_2^\C)$.
\item Set $\xi(t, \l)= g_+(t,\l)^{-1}\xi_0(\l) g_+(t,\l)$.  Then
$$\xi(t,\l)= v_0(t)+ v_{-1}(t)\l^{-1} + \bar v_{-1}(t)\l$$
for  some $v_0(t)\in \fh_0\cap (i\fg_2)$ and $v_{-1}(t)\in
\fh_{-1}$.
\item Set $u_0= -i v_0$, $u_{-1}= -i v_{-1}$, and $k(t)= g_+(t, 1)$.
Then $k(t)\in G_2$ and $u_0, u_{-1}, k$ satisfy \eqref{at} and
\eqref{au}.
   \een
   Moreover, $f(s,t)= e^{As} k_1(t)$ is almost complex in $S^6$, where $k_1(t)$ is the first column of $k(t)$. 
   \ethm

\end{document}